\documentclass{article}
\usepackage{amsmath,amssymb,dsfont,amsfonts}
\usepackage{graphicx,epstopdf}
\usepackage[T1]{fontenc}
\usepackage[utf8]{inputenc}
\usepackage{hyperref}
\def\u{{\mathbf{u}}}
\def\g{{\mathbf{g}}}
\def\n{{\mathbf{n}}}
\def\vv{{\mathbf{v}}}
\def\w{{\mathbf{w}}}
\def\t{{\boldsymbol{\tau}}}
\def\xii{{\boldsymbol{\xi}}}
\numberwithin{equation}{section}
\newtheorem{theorem}{Theorem}[section]
\newtheorem{lemma}{Lemma}[section]

\title{On the Weak Solutions to Steady-State Mixed Navier-Stokes/Darcy Model
\thanks{Subsidized by NSFC(Grant No. 11571274 \& 11171269) and the Ph.D. Programs Foundation of Ministry of Education of China (Grant No. 20110201110027).}}
\author{Yanren Hou\thanks{School of Mathematics and Statistics, Xi'an Jiaotong University,
        Xi'an, Shaanxi 710049, China. ({\tt yrhou@mail.xjtu.edu.cn})} and
    Haibiao Zheng\thanks{Department of Mathematics, East China Normal University,
    	Shanghai Key Laboratory of Pure Mathematics and Mathematical Practice, Shanghai, China. ({\tt hbzheng@math.ecnu.edu.cn})}
}
\date{}
\begin{document}
\maketitle
\setlength{\parindent}{0.5cm}

\begin{abstract}
In this paper, an a priori estimate of weak solutions to the mixed Navier-Stokes/Darcy model with Beavers-Joseph-Saffman's interface condition and the existence of a weak solution are established without the small data and/or the large viscosity restriction for the first time. Based on these results, the global uniqueness of the weak solution is obtained. 
\end{abstract}



\section{Introduction}
Because of the important applications in real world, the mixed Stokes/Darcy and Navier-Stokes/Darcy model received much attention in both theoretical and numerical aspect in last decades. In numerical point of view, coupled finite element methods \cite{Arbo, BADEA,FEA,Chen,Karp,Rui,Urqu}, discontinuous Galerkin methods \cite{CHID,GIRAULT,Kans,Bivi1,Bivi2}, domain decomposition methods \cite{RRDD,Disc1,BB1,BB2,BB3,Disc2,TSDBJ,BDDM,HEX}, Lagrange multiplier methods \cite{LM1,LM2}, interface relaxation methods \cite{IR1,IR2}, and decoupled methods based on two-grid or multi-grid finite element \cite{CAIMX, CAIMU,TGSD,HOU,ZUO0,STG, ZUO1} are extensively studied in last decades. Although there are so many literatures that made great contribution to the numerical simulation of the steady-state mixed Stokes/Darcy and Navier-Stokes/Darcy model with different interface conditions, the existence of a weak solution to the mixed Navier-Stokes/Darcy model with Beavers-Joseph (BJ) interface condition or even with the more simpler Beavers-Joseph-Saffman (BJS) interface condition for general data keeps unresolved. As is pointed out in \cite{GIRAULT}, the difficulty for obtaining an a priori estimate of the weak solutions and deriving the existence of a weak solution without the restriction of the small data and/or the large viscosity, comes from the interface conditions, which  does not completely compensate the nonlinear convection in the energy balance in the Navier-Stokes equations. Hence,  deriving an a priori estimate for solutions of the weak formulation of the Navier-Stokes/Darcy model even with BJS interface condition is still an unsolved question. Therefore, as far as we know, the global uniqueness of the weak solution still remains an unsolved open question. In \cite{GIRAULT}, the authors derived an existence result of a weak solution when the kinematic viscosity of the fluid flow is large and/or the data is small. Then they got a local uniqueness result. Similar results can be found also in \cite{Disc1} and  \cite{CHID}. Some other authors discussed such existence problem by other different manners. For example, the authors of \cite{BADEA} also proved the existence and uniqueness of the weak solution of the model with BJS interface condition in a closed convex subset by means of Steklov-Poincar\'{e} operator   when the data is small; the authors of \cite{HEX} obtained the well-posedness of the model with BJ interface condition in the sense of the branch of the nonsingular solutions, which demands the Fr\'{e}chet derivative of the Navier-Stokes/Darcy operator to be nonsingular and of course will bring some small data and/or large viscosity restriction.   

In this paper, by expanding the Navier-Stokes/Darcy model with BJS interface condition to a more large coupled system, we resolve the open question raised in \cite{GIRAULT} of the derivation of an a priori estimate of weak solutions. By the same method for obtaining the a priori estimate of the weak solutions, we also obtain the existence of a weak solution without the restriction of the small data and/or the large viscosity. Having the a priori estimate of weak solutions at hand, we get the global uniqueness of the weak solution to the coupled system for the first time. 

The rest of this paper is arranged as follows. In  section 2, we give a brief introduction of the Navier-Stokes/Darcy model with BJS interface conditions and its weak forms. In section 3, an a priori estimate of the weak solutions to the system is obtained by introducing an auxiliary differential system that is subject to the Navier-Stokes/Darcy model. Finally, we give the existence result of a weak solution without the small data and/or the large viscosity restriction by the Galerkin method and obtain the global uniqueness of the weak solution.  

\section{Mixed Navier-Stokes/Darcy model with BJS interface condition}
\begin{figure}[ht]
    \centering
    \includegraphics[width=10cm]{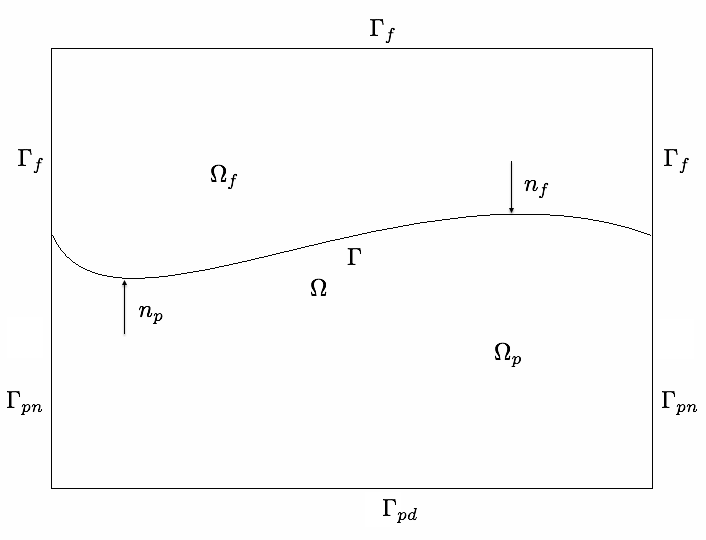}
    \caption{\small A global domain $\Omega$ consisting of a fluid flow region $\Omega_f$ and a porous media flow region $\Omega_p$ separated by an interface $\Gamma$.}\label{fig1}
\end{figure}

Let us consider the following mixed model of the Navier-Stokes equations and the Darcy equation for coupling a fluid flow and a porous media flow in a bounded smooth domain $\Omega\subset {\bf R}^d$, $d=2,3$. Here $\Omega=\Omega_f\cup\Gamma\cup\Omega_p$, where $\Omega_f$ and $\Omega_p$ are two disjoint, connected and bounded domains occupied by fluid flow and porous media flow and $\Gamma=\overline{\Omega}_f\cap\overline{\Omega}_p$ is the interface. We denote $\Gamma_f=\partial\Omega_f\cap\partial\Omega$, $\Gamma_p=\partial\Omega_p\cap\partial\Omega$ and we also denote by $\n_p$ and $\n_f$ the unit outward normal vectors on $\partial\Omega_p$ and $\partial\Omega_f$, respectively. Furthermore, $\Gamma_p$ consists of two disjoint parts $\Gamma_{pd}$ and $\Gamma_{pn}$. We assume $|\Gamma_f|,\;|\Gamma_{pd}|>0$.
See Figure \ref{fig1} for a sketch.

Let us denote by $(\u_f,p_f)$ the velocity field and the pressure of the fluid flow in $\Omega_f$ and $\phi_p$ the piezometric head in $\Omega_p$. The partial differential equations modeling the fluid flow and the porous media flow are 
\begin{equation}\label{NS-Darcy}
\left\{\begin{array}{ll}
-\nabla\cdot(\mathds{T}_\nu(\u_f,p_f))+\u_f\cdot\nabla \u_f=\g_f, & \mbox{in}\;\Omega_f,\\
\nabla\cdot \u_f=0, & \mbox{in}\;\Omega_f,\\
-\nabla\cdot {\mathds{K}}\nabla\phi_p=g_p, & \mbox{in}\;\Omega_p,
\end{array}\right.
\end{equation}
where 
$$
\mathds{T}_\nu(\u_f,p_f)=-p_f\mathds{I}+2\nu\mathds{D}(\u_f),\quad \mathds{D}(\u_f)=\frac12(\nabla \u_f+\nabla^T\u_f),
$$
are the stress tensor and the deformation rate tensor, $\nu>0$ is the kinetic viscosity and $\mathds{K}$ is the permeability in $\Omega_p$, which is a positive definite symmetric tensor that is allowed to vary in space. The third equation of (\ref{NS-Darcy}) that describes the porous media flow motion is the Darcy's law for the piezometric head $\phi_p$. 

The above equations (\ref{NS-Darcy}) are completed and coupled together by the following boundary conditions:
\begin{equation}\label{boundary}
\u_f=0 \quad \mbox{on}\;\Gamma_f,\quad\mathds{K}\nabla\phi_p\cdot \n_p=0\quad\mbox{on}\;\Gamma_{pn},\quad\phi_p=0\quad\mbox{on}\;\Gamma_{pd},
\end{equation}
and the interface conditions on $\Gamma$:
\begin{equation}\label{BJS}
\left\{\begin{array}{l}
\u_f\cdot \n_f-\mathds{K}\nabla\phi_p\cdot \n_p=0,\\
-[\mathds{T}_\nu(\u_f,p_f)\cdot \n_f]\cdot \n_f=\phi_p,\\
-[\mathds{T}_\nu(\u_f,p_f)\cdot \n_f]\cdot\t_i=G\u_f\cdot\t_i,\quad
i=1,\cdots,d-1.
\end{array}\right.
\end{equation}
Here $G>0$ is a constant depending on the nature of the porous medium and usually determined from experimental data, $\t_i$, $i=1,\cdots,d-1$, are the orthonormal tangential unit vectors along $\Gamma$. The first condition is the mass conservation, the second one is the balance of normal force and the third interface condition means the shear force is proportional to the tangential components of the fluid velocity, which is called the Beavers-Joseph-Saffman's (BJS) interface condition (see \cite{BJ} and \cite{BJS}). For more details of these equations, we refer readers to \cite{GIRAULT} and \cite{TGSD}. And in the rest of this paper, we always use boldface characters to denote vectors or vector valued spaces.  

Let us introduce the following Hilbert spaces
\begin{eqnarray*}
&&{\bf X}_f =\{\vv_f\in H^1(\Omega_f)^d: \vv_f|_{\Gamma_f}=0\},\\
&&X_p = \{\psi_p \in H^1(\Omega_p): \psi_p|_{\Gamma_{pd}}=0\},\\
&&Q_f=L^2_0(\Omega_f)=\{q_f\in L^2(\Omega_f): \int_{\Omega_f}q_f=0\},
\end{eqnarray*}
where $({\bf X}_f, Q_f)$ is the space pair for the velocity and the pressure in the fluid flow region $\Omega_f$ and $X_p$ is the space for the piezometric head in the porous medium region $\Omega_p$.

Furthermore, we assume 
\begin{equation}\label{A1}
\g_f\in {\bf X}'_f,\quad g_p\in X'_p,\quad \mathds{K}\in L^\infty(\Omega_p)^{d\times d}
\end{equation}
and there exist $\lambda_{max}>0$ and $\lambda_{min}>0$ such that
\begin{equation}\label{A2}
\mbox{a.e.}\;\mathbf{x}\in \Omega_p,\quad \lambda_{min}|\mathbf{x}|^2\leq \mathds{K}\mathbf{x}\cdot \mathbf{x}\leq\lambda_{max}|\mathbf{x}|^2.
\end{equation}
Here ${\bf X}'_f$ and $X'_p$ are the dual spaces of ${\bf X}_f$ and $X_p$, respectively.

For simplicity, we always use $(\cdot,\cdot)_D$ and $\|\cdot\|_D$ to denote the $L^2$ inner product and the corresponding norm on any given domain $D$. Since $|\Gamma_f|,\;|\Gamma_{pd}|>0$,
we know that $\|D(\cdot)\|_{\Omega_f}$ and $\|\mathds{K}^\frac12\nabla\cdot\|_{\Omega_p}$ are the equivalent norms of the usual Sobolev norms in ${\bf X}_f$ and $X_p$ due to the Korn's inequality, the Poincar\'{e} inequality and (\ref{A2}). 

In addition, let us denote 
$${\bf Y}={\bf X}_f\times Q_f\times X_p.$$ 
Now the weak formulation of the mixed Navier-Stokes/Darcy model with BJS interface condition reads as follows (see \cite{CAIMX}, \cite{GIRAULT}, \cite{LM2} and \cite{TGSD} for details):
for $\g_f\in {\bf X}_f',\;g_p\in X_p'$, find $(\u_f,p_f,\phi_p) \in {\bf Y}$ such that $\forall (\vv_f,q_f,\psi_p) \in{\bf Y}$
\begin{equation*}\left\{\begin{array}{l}
2\nu(D(\u_f),D(\vv_f))_{\Omega_f}+((\u_f\cdot\nabla)\u_f,\vv_f)_{\Omega_f}+\frac12(\nabla\cdot\u_f,\u_f\cdot\vv_f)_{\Omega_f}\\
\qquad+(\mathds{K}\nabla\phi_p,\nabla\psi_p)_{\Omega_p}-(p_f,\nabla\cdot \vv_f)_{\Omega_f}+(q_f,\nabla\cdot \u_f)_{\Omega_f}\\
\qquad+(\phi_p,\vv_f\cdot \n_f)_\Gamma-(\psi_p,\u_f\cdot \n_f)_\Gamma+\sum\limits_{i=1}^{d-1}G(\u_f\cdot\t_i,\vv_f\cdot\t_i)_\Gamma\\
\quad=(\g_f,\vv_f)_{\Omega_f}+(g_p,\psi_p)_{\Omega_p}.
\end{array}\right.\leqno{(Q)}
\end{equation*}
Since $\nabla\cdot \u_f=0$, we know that $\frac12(\nabla\cdot \u_f, \u_f\cdot\vv_f)_{\Omega_f}=0$ in the problem ($Q$).

Thanks to \cite{GIRAULT}, we know that there exists a positive constant $\beta>0$ such that the following Ladyzhenskaya-Babu\v{s}ka-Brezzi (LBB) condition holds:
\begin{equation}\label{LBB}
\inf\limits_{q_f\in Q_f}\sup\limits_{\vv_{f}\in {\bf X}_f}\frac{(q_f,\nabla\cdot \vv_f)_{\Omega_f}}{\|q_f\|_{Q_f}\,\|\vv_f\|_{{\bf X}_f}}\geq \beta.
\end{equation}

If we introduce the following divergence-free space
$$
{\bf V}_f=\{\vv_f\in {\bf X}_f:\;\nabla\cdot \vv_f=0\},
$$
and 
$$
{\bf W}={\bf V}_f\times X_p,
$$
the  restriction of the test function $\vv_f$ to ${\bf V}_f$ in ($Q$) leads to the following reduced weak form: find $(\u_f,\phi_p)\in {\bf W}$ such that $\forall (\vv_f,\psi_p)\in{\bf W}$
\begin{equation*}\left\{\begin{array}{l}
2\nu(D(\u_f),D(\vv_f))_{\Omega_f}+((\u_f\cdot\nabla)\u_f,\vv_f)_{\Omega_f}+\frac12(\nabla\cdot \u_f, \u_f\cdot\vv_f)_{\Omega_f}\\
\qquad+(\mathds{K}\nabla\phi_p,\nabla\psi_p)_{\Omega_p}+(\phi_p,\vv_f\cdot \n_f)_\Gamma-(\psi_p,\u_f\cdot \n_f)_\Gamma\\
\qquad+\sum\limits_{i=1}^{d-1}G(\u_f\cdot\t_i,\vv_f\cdot\t_i)_\Gamma=(\g_f,\vv_f)_{\Omega_f}+(g_p,\psi_p)_{\Omega_p}.
\end{array}\right.\leqno{(P)}
\end{equation*}
By the same argument in \cite{GIRAULTBOOK}, we know that the problem ($Q$) and ($P$) are equivalent. 

For the purpose of later analysis, we recall some inequalities and identity:
 \begin{eqnarray}
 \label{L2G}&&\|v\|_{L^2(\partial D)}\leq c\|v\|^\frac12_{L^2(D)}\|v\|^\frac12_{H^1(D)}\leq c\|v\|_{H^1(D)}\quad\forall v\in H^1(D),\\
 \label{L4G}&&\|v\|_{L^4(\partial D)}\leq c\|v\|_{H^1(D)}\quad\forall v\in H^1(D),\\
 \label{tri}&& |((\w_f\cdot\nabla)\u_f,\vv_f)_{\Omega_f}|,\;|(\nabla\cdot \u_f,\w_f\cdot \vv_f)_{\Omega_f}|\\
  \nonumber &&\qquad\leq c\|\w_f\|_{\Omega_f}^\frac12\|D(\w_f)\|_{\Omega_f}^\frac12\|D(\u_f)\|_{\Omega_f}\|D(\vv_f)\|_{\Omega_f},\quad\forall \w_f,\u_f,\vv_f\in {\bf X}_f,\\
\label{antisy}&& \int_D(\u\cdot\nabla)\vv\cdot \w+\int_D(\u\cdot\nabla)\w\cdot \vv\\
\nonumber &&\qquad =\int_{\partial D}(\vv\cdot \w)\u\cdot \n-\int_D(\vv\cdot \w)\nabla\cdot \u, \quad\forall \u,\vv,\w\in H^1(D)^d.
\end{eqnarray}
Here and after, we always use $c$ to denote a generic positive constant which may take different values in different occasions.

\section{An a priori estimate of weak solutions}
In this section, we will establish an a priori estimate for solutions to ($P$).
From now on, we assume that $\Omega_f$, $\Omega_p$ and $\Omega$ are all polygons or polyhedrons for simplicity. For a given small positive parameter $h>0$, let us denote by $T^h_f$, $T^h_p$ the regular triangulations of $\Omega_f$, $\Omega_p$ and we assume that the two meshes coincide on $\Gamma$. It is obvious that the mergence of the above two triangulations forms a regular triangulation $T^h$ of $\Omega$. 
We also denote
$$
\Omega_{fh}=\bigcup\limits_{K\in T^h_f} K, \quad \Omega_{ph}=\bigcup\limits_{K\in T^h_p} K,\quad \Omega_h=\bigcup\limits_{K\in T^h} K.
$$
For $h>0$ small enough, we assume that $\Omega_{fh}=\Omega_f$, $\Omega_{ph}=\Omega_p$ and $\Omega_h=\Omega$ and we will not distinguish between $\Omega$, $\Omega_f$, $\Omega_p$ and $\Omega_h$, $\Omega_{fh}$, $\Omega_{ph}$ in the later analysis.

Let us denote by ${\bf X}_h\subset H^1_0(\Omega)^d$ a finite element space defined on $\Omega$ and we denote by ${\bf X}_{fh}={\bf X}_h|_{\Omega_{fh}}\subset{\bf X}_f$ a finite element space defined on $\Omega_f$. And we denote by $Q_{fh}\subset Q_f$ and $X_{ph}\subset X_p$ the finite element spaces based on the above triangulations and we assume that $({\bf X}_{fh}, Q_{fh})$ is a stable finite element pair. 
Furthermore, we define the following vector valued Hilbert space on $\Omega_p$
$$
{\bf X}_p=\{\vv_p\in H^1(\Omega_p)^d: \vv_p|_{\Gamma_p}=0\}.
$$

Now we introduce a finite element space ${\bf V}_h\subset{\bf X}_h$ defined on $\Omega$ as
$$
{\bf V}_h=\{\vv_h\in {\bf X}_h: (\nabla\cdot \vv_h, q_{fh})_{\Omega_f}=0,\;\forall q_{fh}\in Q_{fh}\}.
$$
And we define
$$
{\bf V}_{fh}={\bf V}_h|_{\Omega_f}\subset{\bf X}_f,\quad {\bf X}_{ph}={\bf V}_h|_{\Omega_p}\subset{\bf X}_p,
$$
where ${\bf V}_{fh}$ is a weakly divergence-free finite element space defined on $\Omega_f$ and ${\bf X}_{ph}$ is a vector valued finite element space defined on $\Omega_p$. In the rest, we also denote
$$
{\bf W}_h={\bf V}_{fh}\times X_{ph}.
$$

In addition, let us denote by $\Pi_f^h$ the Scott-Zhang interpolator\cite{SCOTTZHANG} from ${\bf X}_f$ onto ${\bf X}_{fh}$ with the following property
\begin{equation}\label{interpolation}
\|\vv_f-\Pi_f^h\vv_f\|_{\Omega_f}\leq ch\|D(\vv_f)\|_{\Omega_f},\quad\forall \vv_f\in {\bf X}_f.
\end{equation}

We know from \cite{GIRAULT} that the difficulty for obtaining an a priori estimate of the Navier-Stokes/Darcy model comes from the unbalance of the energy cause by the nonlinear convection in the Navier-Stokes equations. This observation suggests us to construct an auxiliary partial differential system defined in $\Omega_p$, which is subjected to the Navier-Stokes/Darcy problem, so that the auxiliary system can completely or almost completely compensate the nonlinear convection in the energy balance of the Navier-Stokes equations. To do this, we first introduce a lifting operator $\gamma^{-1}$ from $(H_{00}^\frac12(\Gamma))^d$ to ${\bf X}_p$ (see \cite{GIRAULTBOOK}): for any $\boldsymbol{\zeta}\in (H_{00}^\frac12(\Gamma))^d$ with 
$\int_{\Gamma}\boldsymbol{\zeta}ds=0$, 
$$
\gamma^{-1}\boldsymbol{\zeta}\in{\bf X}_p,\quad (\gamma^{-1}\boldsymbol{\zeta})|_\Gamma=\boldsymbol{\zeta},\quad \nabla\cdot(\gamma^{-1}\boldsymbol{\zeta})=0.
$$
Now we define the following auxiliary linear equations in $\Omega_p$: for any given $\xii\in {\bf V}_f$
\begin{equation}\label{auxe}
\left\{\begin{array}{l}
-2\sigma\nabla\cdot \mathds{D}(\u_p)+(\u_p^0\cdot\nabla)\u_p=0\quad\mbox{in}\;\Omega_p,\\
\u_p|_{\Gamma_p}=0,\quad \u_p|_{\Gamma}=\xii_f|_\Gamma.
\end{array}\right.
\end{equation}
where $\u_p^0=\gamma^{-1}(\xii_f|_\Gamma) \in {\bf X}_p$
and $\sigma>0$ is a certain positive constant which will be specified later.
It is obvious that this linear auxiliary system for any given $\sigma>0$ and $\xii_f\in {\bf V}_f$ is well-posed. 

Now, for $\xii_f=\u_f\in {\bf V}_f$, we consider its Galerkin approximation in the finite element space ${\bf X}_{ph}$, which is:
for $\u_f\in {\bf V}_f$, find $\u_{ph}\in {\bf X}_{ph}$ such that $\forall \vv_{ph}\in {\bf X}_{ph}$
\begin{equation}\label{auxns}
\left\{\begin{array}{l}
2\sigma(D(\u_{ph}),D(\vv_{ph}))_{\Omega_p}+((\u_p^0\cdot\nabla)\u_{ph},\vv_{ph})_{\Omega_p}
-\sigma\int_\Gamma\frac{\partial \u_{ph}}{\partial \n_p}\vv_{ph}=0,\\
\u_{ph}|_\Gamma=\Pi_f^h\u_f|_\Gamma,
\end{array}\right.
\end{equation}
where $\u_p^0=\gamma^{-1}(\u_f|_\Gamma)$. 
If we introduce the following finite element space
$$
{\bf \mathring{X}}_{ph}=\{\vv_{ph}\in{\bf X}_{ph}: \vv_{ph}|_\Gamma=0\}\subset H^1_0(\Omega_p)^d,
$$
an equivalent form of (\ref{auxns}) is: find $\u_{ph}\in {\bf X}_{ph}$ such that $\forall \vv_{ph}\in {\bf \mathring{X}}_{ph}$
\begin{equation}\label{auxns2}
\left\{\begin{array}{l}
2\sigma(D(\u_{ph}),D(\vv_{ph}))_{\Omega_p}+((\u_p^0\cdot\nabla)\u_{ph},\vv_{ph})_{\Omega_p}=0,\\
\u_{ph}|_\Gamma=\Pi_f^h\u_f|_\Gamma,
\end{array}\right.
\end{equation}

Now it is ready for us to derive the a priori estimate of solutions to ($P$).

\begin{theorem}\label{th0} There holds the following a priori estimate for solutions $(\u_f,\phi_p)\in {\bf W}$ to ($P$) 
$$\nu\|D(\u_f)\|_{\Omega_f}^2+\|\mathds{K}^\frac12\nabla\phi_p\|_{\Omega_p}^2\leq{\cal C}^2,$$
where 
$$
{\cal C}^2=c\nu^{-1}\|\g_f\|_{{\bf X}'_f}^2+c\lambda_{min}^{-1}\|g_p\|_{X_p'}^2,
$$
and $c>0$ is a generic constant that has nothing to do with the data of the problem.
\end{theorem}
\noindent{\bf Proof.}  It is clear that ($P$) and (\ref{auxns}) form a new coupled system, in which (\ref{auxns}) is subjected to ($P$) while ($P$) has nothing to do with (\ref{auxns}). 

In the rest of the proof, we assume that $(\u_f,\phi_p)\in {\bf W}$ is a possible solution to the problem ($P$) and there exists a positive constant $M_{\u_f}<\infty$ such that 
\begin{equation}\label{A3}
\|D(\u_f)\|_{\Omega_f}\leq M_{\u_f}.
\end{equation}

Taking $(\vv_f,\psi_p)=(\u_f,\phi_p)$ in ($P$) and omitting the non-negative term $\sum\limits_{i=1}^{d-1}G(\u_f\cdot\t_i,\u_f\cdot\t_i)_\Gamma$, we have
\begin{eqnarray}\label{GAE1}
&&2\nu\|D(\u_f)\|_{\Omega_f}^2+\|\mathds{K}^\frac12\nabla\phi_p\|_{\Omega_p}^2+((\u_f\cdot\nabla)\u_f,\u_f)_{\Omega_f}\\
&&\nonumber\qquad+\frac12(\nabla\cdot \u_f, \u_f\cdot\u_f)_{\Omega_f}\leq(\g_f,\u_f)_{\Omega_f}+(g_p,\phi_p)_{\Omega_p}.
\end{eqnarray}
Taking $\vv_{ph}=\u_{ph}$ in (\ref{auxns}) and using the boundary condition $\u_{ph}|_\Gamma=\Pi_f^h\u_f|_\Gamma$ lead to
\begin{eqnarray}\label{GAE2}
&&2\sigma\|D(\u_{ph})\|_{\Omega_p}^2+((\u_p^0\cdot\nabla)\u_{ph},\u_{ph})_{\Omega_p}-\sigma\int_\Gamma\frac{\partial \u_{ph}}{\partial \n_p}\Pi_f^h\u_f=0.
\end{eqnarray}
Being aware of $\n_f=-\n_p$ on $\Gamma$, $\nabla\cdot\u_f=\nabla\cdot\u_p^0=0$ and the identity (\ref{antisy}), it is easy to verify that
\begin{eqnarray*}
&&((\u_f\cdot\nabla)\u_f,\u_f)_{\Omega_f}+\frac12(\nabla\cdot \u_f, \u_f\cdot\u_f)_{\Omega_f}=\frac12\int_\Gamma|\u_f|^2\u_f\cdot \n_f,\\
&&((\u_p^0\cdot\nabla)\u_{ph},\u_{ph})_{\Omega_p}=\frac12\int_\Gamma|\u_{ph}|^2\u_p^0\cdot \n_p=-\frac12\int_\Gamma|\Pi_f^h\u_f|^2\u_f\cdot \n_f.
\end{eqnarray*}

By using (\ref{L2G}), (\ref{L4G}), (\ref{interpolation}), the assumption (\ref{A3}), the Korn's and the Poincar\'{e} inequality, summation of the above two identities leads to
\begin{eqnarray}
\label{nonconv}&&((\u_f\cdot\nabla)\u_f,\u_f)_{\Omega_f}+\frac12(\nabla\cdot \u_f, \u_f\cdot\u_f)_{\Omega_f}+((\u_p^2\cdot\nabla)\u_{ph},\u_{ph})_{\Omega_p}\\
\nonumber&&\quad=\frac12\int_\Gamma[|\u_f|^2\u_f\cdot \n_f-|\Pi_f^h\u_f|^2\u_f\cdot \n_f]\\
\nonumber&&\quad=\frac12\int_\Gamma[(\u_f-\Pi_f^h\u_f)\cdot(\u_f+\Pi_f^h\u_f)\u_f\cdot\n_f]\\
\nonumber&&\quad\leq c\|\u_f+\Pi_f^h\u_f\|_{L^4(\Gamma)}\|\u_f\|_{L^4(\Gamma)})\|\u_f-\Pi_f^h\u_f\|_{L^2(\Gamma)}\\
\nonumber&&\quad\leq cM_{\u_f}h^\frac12\|D(\u_f)\|_{\Omega_f}^2.
\end{eqnarray}
This means the auxiliary system can almost compensate the the nonlinear convection of the Navier-Stokes equations in the energy balance as $h\rightarrow 0$.

Taking the above estimation into account, the summation of (\ref{GAE1}) and (\ref{GAE2}) yields
\begin{eqnarray}\label{INEQ1}
&&2\nu\|D(\u_f)\|_{\Omega_f}^2+\|\mathds{K}^\frac12\nabla\phi_p\|_{\Omega_p}^2+2\sigma\|D(\u_{ph})\|_{\Omega_p}^2\\
&&\nonumber\qquad\leq |(\g_f,\u_f)_{\Omega_f}|+|(g_p,\phi_p)_{\Omega_p}|+\sigma|\int_\Gamma\frac{\partial \u_{ph}}{\partial \n_p}\Pi_f^h\u_f|+cM_{\u_f}h^\frac12\|D(\u_f)\|_{\Omega_f}^2.
\end{eqnarray}

For the first and the second term on the right hand side of the above inequality (\ref{INEQ1}), by using the Korn's inequality, the Poincar\'{e} inequality and (\ref{A2}) we have
\begin{eqnarray}
\label{right2}&&|(\g_f,\u_f)_{\Omega_f}|+|(g_p,\phi_p)_{\Omega_p}|\\
\nonumber&&\qquad\leq c\|\g_f\|_{{\bf X}'_f}\|D(\u_f)\|_{\Omega_f}+c\lambda_{min}^{-\frac12}\|g_p\|_{X'_p}\|\mathds{K}^\frac12\nabla\phi_p\|_{\Omega_p}\\
\nonumber&&\qquad\leq\frac{\nu}{2}\|D(\u_f)\|_{\Omega_f}^2+\frac12\|\mathds{K}^\frac12\nabla\phi_p\|_{\Omega_p}^2+c\nu^{-1}\|\g_f\|_{{\bf X}'_f}^2+c\lambda_{min}^{-1}\|g_p\|_{X'_p}^2.
\end{eqnarray}

For the third term on the right hand side of (\ref{INEQ1}), by using (\ref{L2G}), the Korn's inequality, the Poincar\'{e} inequality and the following inequality (see \cite{CIARLET} and \cite{PAUL})
$$
\|v\|_{L^2(\partial K)}\leq ch^{-\frac12}\|v\|_{L^2(K)},\quad\mbox{for any polynomial}\;v\;\mbox{on}\; K,
$$
we have
\begin{eqnarray}
\label{bint}&&\sigma|\int_{\Gamma}\frac{\partial \u_{ph}}{\partial \n_p}\Pi_f^h\u_f|\leq\sigma\|\frac{\partial \u_{ph}}{\partial \n_p}\|_{L^2(\Gamma)}\|\Pi_f^h\u_f\|_{L^2(\Gamma)}\\
\nonumber&&\qquad\leq c\sigma\|\nabla\u_{ph}\cdot \n_p\|_{L^2(\partial\Omega_p)}\|\Pi_f^h\u_f\|_{L^2(\partial\Omega_f)}\\
\nonumber&&\qquad\leq c\sigma(\sum\limits_{K\in T^h_p}\|\nabla\u_{ph}\|^2_{L^2(\partial K)})^\frac12\|\Pi_f^h\u_f\|_{H^1(\Omega_f)}\\
\nonumber&&\qquad\leq c\sigma h^{-\frac12}\|\nabla \u_{ph}\|_{\Omega_p}\|\nabla \u_f\|_{\Omega_f}\leq c\sigma h^{-\frac12}\|D(\u_{ph})\|_{\Omega_p}\|D(\u_{f})\|_{\Omega_f}\\
\nonumber&&\qquad\leq c\sigma^2 h^{-1}\nu^{-1}\|D(\u_{ph})\|_{\Omega_p}^2+\frac{\nu}{2}\|D(\u_{f})\|_{\Omega_f}^2,
\end{eqnarray}
where $\n_K$ stands for the unit outward vector of each element $K\in T^h_p$.

If we choose $h$ small enough and $\sigma$ small enough such that 
$$cM_{\u_f}h^\frac12<\frac{\nu}{2},\quad 0<\sigma\leq c\nu h,$$ 
combination of (\ref{GAE1}), (\ref{GAE2}), (\ref{nonconv}), (\ref{right2}) and (\ref{bint}) admits
\begin{eqnarray}\label{bound}
&&\sigma\|D(\u_{ph})\|_{\Omega_p}^2+\nu\|D(\u_f)\|_{\Omega_f}^2+\|\mathds{K}^\frac12\nabla\phi_p\|_{\Omega_p}^2\\
&&\nonumber\qquad\leq c\nu^{-1}\|\g_f\|_{{\bf X}'_f}^2+c\lambda_{min}^{-1}\|g_p\|_{X'_p}^2\stackrel{\triangle}{=}{\cal C}^2.
\end{eqnarray}

Since the solutions of ($P$) is independent of the system (\ref{auxns}), the above a priori estimate actually gives an $h$ and $\sigma$ independent a priori estimate of solutions to ($P$)
$$\nu\|D(\u_f)\|_{\Omega_f}^2+\|\mathds{K}^\frac12\nabla\phi_p\|_{\Omega_p}^2\leq{\cal C}^2.$$
\hfill$\Box$

\section{Existence and global uniqueness of the weak solution}
In this section, we will use the Galerkin method to show that there exists at least one solution to ($P$) (and ($Q$)), and then give the global uniqueness of the weak solution based on the a priori estimate of weak solutions obtained in last section. 

We first recall the following Brouwer's fixed point theorem.

\noindent{\bf Brouwer's Fixed Point Theorem} {\it Suppose $H$ is a finite dimensional Hilbert space equipped with inner product $(\cdot,\cdot)$ and norm $|\cdot|$. If $\Phi$ is a continuous map from $H$ to $H$ and for a certain constant $\mu>0$
$$
(\Phi(f),f)\geq 0\quad\forall f\in H,\quad |f|=\mu,
$$
there exists a function $f\in H$ such that
$$
\Phi(f)=0,\quad |f|\leq\mu.
$$}

Now let us give the Galerkin approximation of ($P$) in ${\bf W}_h$, the finite element space defined in the previous section: find $(\u_{fh},\phi_{ph})\in {\bf W}_h$ such that $\forall (\vv_{fh},\psi_{ph})\in {\bf W}_h$
\begin{equation*}\left\{\begin{array}{l}
2\nu(D(\u_{fh}),D(\vv_{fh}))_{\Omega_f}+((\u_{fh}\cdot\nabla)\u_{fh},\vv_{fh})_{\Omega_f}\\
\qquad+\frac12(\nabla\cdot\u_{fh},\u_{fh}\cdot\vv_{fh})_{\Omega_f}+(\mathds{K}\nabla\phi_{ph},\nabla\psi_{ph})_{\Omega_p}\\
\qquad+(\phi_{ph},\vv_{fh}\cdot \n_f)_\Gamma-(\psi_{ph},\u_{fh}\cdot \n_f)_\Gamma\\
\qquad+\sum\limits_{i=1}^{d-1}G(\u_{fh}\cdot\t_i,\vv_{fh}\cdot\t_i)_\Gamma=(\g_f,\vv_{fh})_{\Omega_f}+(g_p,\psi_{ph})_{\Omega_p}.
\end{array}\right.\leqno{(P_h)}
\end{equation*}

\begin{lemma}\label{lem0} For any given $h>0$, we have the following a priori estimate for the solutions $(\u_{fh},\phi_{ph})\in {\bf W}_h$ of ($P_h$)
$$\nu\|D(\u_{fh})\|_{\Omega_f}^2+\|\mathds{K}^\frac12\nabla\phi_{ph}\|_{\Omega_p}^2\leq{\cal C}^2,$$
where ${\cal C}$ is defined in Theorem \ref{th0}.
\end{lemma} 

\noindent{\bf Proof.} The proof of this lemma is completely the same as that of Theorem \ref{th0}. In fact, we consider the coupled system ($P_h$) and (\ref{auxns}) with $\xii_f=\u_{fh}$, $\u^0_p=\gamma^{-1}(\u_{fh}|_\Gamma)$ with $\nabla\cdot\u^0_p=0$ this time and we can get similar estimates of (\ref{GAE1}), (\ref{GAE2}), (\ref{right2}) and (\ref{bint}) by just replacing $(\u_f,\phi_p)$ with $(\u_{fh},\phi_{ph})$. For the estimation (\ref{nonconv}), noticing the construction of $\u^0_p$ in (\ref{auxns}) and the identity (\ref{antisy}), it is easy to show
\begin{eqnarray*}
&&((\u_{fh}\cdot\nabla)\u_{fh},\u_{fh})_{\Omega_f}+\frac12(\nabla\cdot\u_{fh},\u_{fh}\cdot\u_{fh})_{\Omega_f}+((\u^0_p\cdot\nabla)\u_{ph},\u_{ph})_{\Omega_p}=0.
\end{eqnarray*}
That is the convection in the auxiliary system (\ref{auxns}) completely compensates the nonlinear convection in the Galerkin approximation of the Navier-Stokes equations. This is also the motivation for us to introduce the auxiliary system (\ref{auxe}). Finally we can get
\begin{eqnarray}\label{boundh}
&&\sigma\|D(\u_{ph})\|_{\Omega_p}^2+\nu\|D(\u_{fh})\|_{\Omega_f}^2+\|\mathds{K}^\frac12\nabla\phi_{ph}\|_{\Omega_p}^2\leq{\cal C}^2,
\end{eqnarray}
for $\sigma$ small enough such that
$$
0<\sigma\leq c\nu h.
$$
This concludes the proof of this lemma.
\hfill$\Box$

Having the a priori estimate (\ref{boundh}) of the coupled system ($P_h$) and (\ref{auxns}), it is ready for us to show the existence of a solution to ($P_h$).

\begin{lemma}\label{lem1} For any given $h>0$, there exists at least one solution $(\u_{fh},\phi_{ph})\in {\bf W}_h$ to ($P_h$). 
\end{lemma} 
\noindent{\bf Proof.} We actually intend to prove this lemma by showing the existence of at least one solution to the coupled system ($P_h$) and (\ref{auxns}) for any $h>0$ and the corresponding sufficiently small $\sigma>0$. Once it is done, we can conclude the existence property of ($P_h$) since (\ref{auxns}) is subjected to ($P_h$) but ($P_h$) has nothing to do with (\ref{auxns}).

First of all, we introduce the following space
$$
{\bf U}_h={\bf V}_h\times X_{ph},
$$
where ${\bf V}_h$ and $X_{ph}$ have already been defined in the previous section.
For any $\vv_h\in {\bf V}_h$, let us denote $\vv_{fh}=\vv_h|_{\Omega_f}\in {\bf V}_{fh}$ and $\vv_{ph}=\vv_h|_{\Omega_p}\in {\bf X}_{ph}$. 
We introduce a mapping: ${\cal F}_h:{\bf U}_h\rightarrow {\bf U}_h$, defined for all $(\vv_h,\psi_{ph})\in {\bf U}_h$ by the following
\begin{eqnarray*}
&&\forall (\w_h,\chi_{ph})\in {\bf U}_h,\quad ({\cal F}_h((\vv_h,\psi_{ph}),(\w_h,\chi_{ph}))_{{\bf U}_h}=2\nu(D(\vv_{fh}),D(\w_{fh}))_{\Omega_f}\\
&&\quad+(\mathds{K}\nabla\psi_{ph},\nabla\chi_{ph})_{\Omega_p}+\sigma(D(\vv_{ph}),D(\w_{ph}))_{\Omega_p}+((\tilde\vv_h\cdot\nabla)\vv_h,\w_h)_\Omega\\
&&\quad+\frac12(\nabla\cdot\tilde\vv_h,\vv_h\cdot\w_h)_\Omega+(\psi_{ph},\w_{fh}\cdot \n_f)_\Gamma-(\chi_{ph},\vv_{fh}\cdot \n_f)_\Gamma\\
&&\quad-\sigma\int_\Gamma\frac{\partial \vv_{ph}}{\partial \n_p}\w_{ph}+\sum\limits_{i=1}^{d-1}G(\vv_{fh}\cdot\t_i,\w_{fh}\cdot\t_i)_\Gamma-(\g_f,\w_{fh})_{\Omega_f}-(g_p,\chi_{ph})_{\Omega_p}.
\end{eqnarray*}
Here $\w_{fh}=\w_h|_{\Omega_f}\in {\bf V}_{fh}$ and $\w_{ph}=\w_h|_{\Omega_p}\in {\bf X}_{ph}$. And $\tilde\vv_h$ is defined as 
$$
\tilde\vv_h|_{\Omega_f}=\vv_{fh},\quad\tilde\vv_h|_{\Omega_p}=\gamma^{-1}(\vv_{fh}|_\Gamma)\in{\bf X}_p\quad\mbox{with }\quad\nabla\cdot(\tilde\vv_h|_{\Omega_p})=0.
$$
It is clear that, if $(\vv_h,\psi_{ph})\in{\bf U}_h$ is a zero point of ${\cal F}_h$, we can assert that $\vv_{ph}$ satisfies (\ref{auxns}) with $\vv_{ph}|_\Gamma=\vv_{fh}|_\Gamma$. In fact, if we choose $\chi_{ph}=0$ and $\w_h\in {\bf V}_h$ with $\w_{fh}=0$ and $\w_{ph}\in {\bf \mathring{X}}_{ph}$, we have $\vv_{ph}\in {\bf X}_{ph}$ satisfies (\ref{auxns2}). That is $\vv_{ph}$ satisfies (\ref{auxns}). Therefore we have $(\vv_{fh},\psi_{ph})\in {\bf W}_h$ is a solution of problem ($P_h$). 

Thanks to the previous estimation (\ref{boundh}), we know that for sufficiently small $\sigma>0$
$$
({\cal F}_h((\vv_h,\psi_{ph}),(\vv_h,\psi_{ph}))_{{\bf U}_h}\geq \sigma\|D(\vv_{ph})\|_{\Omega_p}^2+\nu\|D(\vv_{fh})\|_{\Omega_f}^2+\|\mathds{K}^\frac12\nabla\psi_{ph}\|_{\Omega_p}^2-{\cal C}^2.
$$
By equipping the space ${\bf U}_h$ with the following norm
$$
\||(\vv_h,\psi_{ph})\||_\Omega=(\sigma\|D(\vv_{ph})\|_{\Omega_p}^2+\nu\|D(\vv_{fh})\|_{\Omega_f}^2+\|\mathds{K}^\frac12\nabla\psi_{ph}\|_{\Omega_p}^2)^\frac12,
$$
the above estimate ensures 
$$
({\cal F}_h(\vv_h,\psi_{ph}),(\vv_h,\psi_{ph}))_{{\bf U}_h}\geq 0\quad\forall (\vv_h,\psi_{ph})\in {\bf U}_h,\quad \||(\vv_h,\psi_{ph})\||={\cal C}.
$$
Now by using the Brouwer's fixed point theorem, we can conclude the proof of this lemma.\hfill$\Box$

Thanks to the result of Lemma \ref{lem0} and \ref{lem1}, we get a bounded sequence $\{\u_{fh},\phi_{ph}\}_{h>0}$ in ${\bf X}_f\times X_p$. Since ${\bf X}_f\times X_p$ is compactly embeded in $L^2(\Omega_f)^d\times L^2(\Omega_p)$, 
we can extract a subsequence, which is still denoted by $h$, such that as $h\rightarrow 0$ there exists $(\u_f,\phi_p)\in {\bf X}_f\times X_p$ such that
\begin{eqnarray}
\label{weakconv}&&(\u_{fh},\phi_{ph})\longrightarrow (\u_f,\phi_p)\quad\mbox{weakly},\\
\label{strongconv}&& (\u_{fh},\phi_{ph})\longrightarrow (\u_f,\phi_p)\quad\mbox{strongly in}\;L^2(\Omega_f)^d\times L^2(\Omega_p).
\end{eqnarray}

\begin{lemma}\label{lem2} For $(\u_f,\phi_p)$ defined in (\ref{weakconv}) and (\ref{strongconv}), we have
$\forall (\vv_f,\psi_p)\in {\bf W}$
\begin{eqnarray*}
&&\lim\limits_{h\rightarrow 0}\{2\nu(D(\u_{fh}),D(\vv_f))_{\Omega_f}+(\mathds{K}\nabla\phi_{ph},\nabla\psi_p)_{\Omega_p}\}\\
&&\qquad=2\nu(D(\u_f),D(\vv_f))_{\Omega_f}+(\mathds{K}\nabla\phi_p,\nabla\psi_p)_{\Omega_p},\\
&&\lim\limits_{h\rightarrow 0}\{(\phi_{ph},\vv_f\cdot \n_f)_\Gamma-(\psi_p,\u_{fh}\cdot \n_f)_\Gamma+\sum\limits_{i=1}^{d-1}G(\u_{fh}\cdot\t_i,\vv_f\cdot\t_i)_\Gamma\}\\
&&\qquad=(\phi_p,\vv_f\cdot \n_f)_\Gamma-(\psi_p,\u_f\cdot \n_f)_\Gamma+\sum\limits_{i=1}^{d-1}G(\u_f\cdot\t_i,\vv_f\cdot\t_i)_{\Gamma},\\
&&\lim\limits_{h\rightarrow 0}[((\u_{fh}\cdot\nabla)\u_{fh},\vv_f)_{\Omega_f}+\frac12(\nabla\cdot\u_{fh},\u_{fh}\cdot \vv_f)_{\Omega_f}]=((\u_f\cdot\nabla)\u_f,\vv_f)_{\Omega_f}.
\end{eqnarray*}
\end{lemma}

\noindent{\bf Proof.}  For any $(\vv_f,\psi_p)\in {\bf W}$, it is obvious that
\begin{eqnarray*}
&&\lim\limits_{h\rightarrow 0}\{2\nu(D(\u_{fh}),D(\vv_f))_{\Omega_f}+(\mathds{K}\nabla\phi_{ph},\nabla\psi_p)_{\Omega_p}\}\\
&&\qquad=2\nu(D(\u_f),D(\vv_f))_{\Omega_f}+(\mathds{K}\nabla\phi_p,\cdot\nabla\psi_p)_{\Omega_p},
\end{eqnarray*}
because of (\ref{weakconv}). For the second limit, we have
\begin{eqnarray*}
&&\sum\limits_{i=1}^{d-1}G\int_\Gamma(\u_{fh}\cdot\t_i)(\vv_f\cdot\t_i)+\int_\Gamma[\phi_{ph}\vv_f-\psi_p \u_{fh}]\cdot \n_f\\
&&\qquad=\sum\limits_{i=1}^{d-1}G\int_\Gamma((\u_{fh}-\u_f)\cdot\t_i)(\vv_f\cdot\t_i)\\
&&\qquad\quad+\int_\Gamma[(\phi_{ph}-\phi_p)\vv_f-\psi_p(\u_{fh}-\u_f)]\cdot \n_f\\
&&\qquad\quad+\sum\limits_{i=1}^{d-1}G\int_\Gamma(\u_f\cdot\t_i)(\vv_f\cdot\t_i)+\int_\Gamma[\phi_p\vv_f-\psi_p \u_f]\cdot \n_f.
\end{eqnarray*}
For the first and the second term on the right hand side of the above identity, by using (\ref{L2G}), we have the following estimations.
\begin{eqnarray*}
&&|\sum\limits_{i=1}^{d-1}G\int_\Gamma((\u_{fh}-\u_f)\cdot\t_i)(\vv_f\cdot\t_i)|\leq c\|\u_{fh}-\u_f\|_{L^2(\Gamma)}\|\vv_f\|_{L^2(\Gamma)}\\
&&\qquad\leq c\|\u_{fh}-\u_f\|_{L^2(\partial\Omega_f)}\|\vv_f\|_{L^2(\partial\Omega_f)}\\
&&\qquad\leq c\|\u_{fh}-\u_f\|_{\Omega_f}^\frac12\|D(\u_{fh}-\u_f)\|_{\Omega_f}^\frac12\|D(\vv_f)\|_{\Omega_f},\\
&&|\int_\Gamma[(\phi_{ph}-\phi_p)\vv_f-\psi_p(\u_{fh}-\u_f)]\cdot \n_f|\\
&&\qquad\leq c\|\phi_{ph}-\phi_p\|_{L^2(\Gamma)}\|\vv_f\|_{L^2(\Gamma)}+c\|\psi_p\|_{L^2(\Gamma)}\|\u_{fh}-\u_f\|_{L^2(\Gamma)}\\
&&\qquad=c\|\phi_{ph}-\phi_p\|_{L^2(\partial\Omega_p)}\|\vv_f\|_{L^2(\partial\Omega_f)}+c\|\psi_p\|_{L^2(\partial\Omega_p)}\|\u_{fh}-\u_f\|_{L^2(\partial\Omega_f)}\\
&&\qquad\leq c\|\phi_{ph}-\phi_p\|_{\Omega_p}^\frac12\|\nabla(\phi_{ph}-\phi_p)\|_{\Omega_p}^\frac12\|D(\vv_f)\|_{\Omega_f}\\
&&\qquad\quad+\|\nabla\psi_p\|_{\Omega_p}\|\u_{fh}-\u_f\|_{\Omega_f}^\frac12\|D(\u_{fh}-\u_f)\|_{\Omega_f}^\frac12.
\end{eqnarray*}
Thanks to the uniform boundedness of $(\u_{fh},\phi_{ph})$ in $H^1$ norm and (\ref{strongconv}), we know that these two terms tend to zero when $h\rightarrow 0$. Then we obtain
\begin{eqnarray*}
&&\lim\limits_{h\rightarrow 0}\{\sum\limits_{i=1}^{d-1}G\int_\Gamma(\u_{fh}\cdot\t_i)(\vv_f\cdot\t_i)+\int_\Gamma[\phi_{ph}\vv_f-\psi_p \u_{fh}]\cdot \n_f\}\\
&&\qquad=\sum\limits_{i=1}^{d-1}G\int_\Gamma(\u_f\cdot\t_i)(\vv_f\cdot\t_i)+\int_\Gamma[\phi_p\vv_f-\psi_p \u_f]\cdot \n_f.
\end{eqnarray*}

For the third limit, the limit of the trilinear form, by using (\ref{tri}) we have
\begin{eqnarray*}
&&|((\u_{fh}\cdot\nabla)\u_{fh},\vv_f)_{\Omega_f}+\frac12(\u_{fh}\cdot \vv_f,\nabla\cdot \u_{fh})_{\Omega_f}-((\u_f\cdot\nabla)\u_f,\vv_f)_{\Omega_f}|\\
&&\qquad\leq|(((\u_{fh}-\u_f)\cdot\nabla)\u_{fh},\vv_f)_{\Omega_f}|+\frac12|((\u_{fh}-\u_f)\cdot \vv_f,\nabla\cdot \u_{fh})_{\Omega_f}|\\
&&\qquad\quad+|((\u_f\cdot\nabla)(\u_{fh}-\u_f),\vv_f)_{\Omega_f}|+\frac12|(\u_f\cdot \vv_f,\nabla\cdot(\u_{fh}-\u_f))_{\Omega_f}|\\ 
&&\qquad\leq  c\|\u_{fh}-\u_f\|_{\Omega_f}^\frac12\|D(\u_{fh}-\u_f)\|_{\Omega_f}^\frac12\|D(\u_{fh})\|_{\Omega_f}\|D(\vv_f)\|_{\Omega_f}\\
&&\qquad\quad+|((\u_f\cdot\nabla)(\u_{fh}-\u_f),\vv_f)_{\Omega_f}|+|(\u_f\cdot \vv_f,\nabla\cdot(\u_{fh}-\u_f))_{\Omega_f}|.
\end{eqnarray*}
By using the uniform boundedness of $\u_{fh}$ in ${\bf X}_f$ and (\ref{strongconv}) again, the first term on the right hand side of the above last inequality tends to zero as $h$ goes to zero. And the rest two terms on the right hand side of the above last inequality tend to zero when $h\rightarrow 0$ because of (\ref{weakconv}). Then we can derive the third limit in this lemma.
\hfill$\Box$

With the above three lemmas, we could obtain the following existence result of the problem ($P$) and ($Q$).

 \begin{theorem}\label{th1}
The problem ($Q$) has at least one solution $(\u_f,p_f,\phi_p)\in{\bf Y}$ with the following bounds
$$\nu\|D(\u_f)\|_{\Omega_f}^2+\|\mathds{K}^\frac12\nabla\phi_p\|_{\Omega_p}^2\leq{\cal C}^2,$$
and
$$
\|p_f\|_{\Omega_f}\leq c\beta^{-1}{\cal C}(1+{\cal C}),
$$
where ${\cal C}$ is defined in Theorem \ref{th0}.
\end{theorem}

\noindent{\bf Proof.} Taking $h\rightarrow 0$ in problem ($P_h$) and being aware of the results in Lemma \ref{lem2}, we know that the limit $(\u_f,\phi_p)$ in (\ref{weakconv}) and (\ref{strongconv}) is a solution of problem ($P$) and its bound is obvious thanks to Lemma \ref{lem0}.

Thanks to the LBB condition (\ref{LBB}), we know that there exists a unique $p_f\in Q_f$ such that $(\u_f,p_f,\phi_p)\in{\bf Y}$ is a solution of problem ($Q$). And it is easily obtained by using (\ref{LBB}) that
$$
\|p_f\|_{\Omega_f}\leq c\beta^{-1}{\cal C}(1+{\cal C}).
$$
\hfill$\Box$

\noindent{\bf Remark 1} In the above theorem, we get the existence result of a weak solution to the Navier-Stokes/Darcy problem with BJS interface condition in general circumstances. Compared with other existence result, for example the existence result in \cite{GIRAULT}, we get an a priori estimate for weak solutions of problem ($P$) and ($Q$) and the existence property need no small data and/or large viscosity restriction as long as the assumptions (\ref{A1}) and (\ref{A2}) are fulfilled. On the other hand, the a priori estimate for the solutions to problem ($P$) makes it possible to establish the global uniqueness of the solution of the Navier-Stokes/Darcy problem, which is an open problem raised in \cite{GIRAULT}.

\begin{theorem}\label{th2} Assume the data of the mixed Navier-Stokes/Darcy model with BJS interface condition satisfy
\begin{equation}\label{10050959}
c(\nu^{-2}\|\g_f\|_{{\bf X}'_f}+\nu^{-\frac32}\lambda_{min}^{-\frac12}\|g_p\|_{X'_p})<1.
\end{equation}
The problem ($P$) (or ($Q$)) has only one solution in ${\bf W}$ (or ${\bf Y}$).
\end{theorem}

\noindent{\bf Proof.} Suppose $(\u_f^i,\phi_p^i)\in {\bf W}$, $i=1,2$, be two solutions to ($P$). Their difference satisfies
\begin{eqnarray*}
&& 2\nu(D(\u_f^1-\u_f^2),D(\vv_f))_{\Omega_f}+(((\u_f^1-\u_f^2)\cdot\nabla)\u_f^1,\vv_f)_{\Omega_f}\\
&&\qquad+((\u_f^2\cdot\nabla)(\u_f^1-\u_f^2),\vv_f)_{\Omega_f}+(\mathds{K}\nabla(\phi_p^1-\phi_p^2),\nabla\psi_p)_{\Omega_p}\\
&&\qquad+(\phi_p^1-\phi_p^2,\vv_f\cdot \n_f)_\Gamma-(\psi_p,(\u_f^1-\u_f^2)\cdot \n_f)_\Gamma\\
&&\qquad+\sum\limits_{i=1}^{d-1}G\int_\Gamma((\u_f^1-\u_f^2)\cdot\t_i)(\vv_f\cdot\t_i)=0.
\end{eqnarray*}

Taking $\vv_f=\u_f^1-\u_f^2$ and $\psi_p=\phi_p^1-\phi_p^2$, using (\ref{tri}) and taking the result in Theorem \ref{th1} into account , we obtain
\begin{eqnarray*}
&&2\nu\|D(\u^1_f-\u^2_f)\|^2_{\Omega_f}+\|\mathds{K}^\frac12\nabla(\phi_p^1-\phi_p^2)\|^2_{\Omega_p}\\
&&\qquad\leq c(\|D(\u^1_f)\|_{\Omega_f}+\|D(\u^2_f)\|_{\Omega_f})\|D(\u^1_f-\u^2_f)\|^2_{\Omega_f}\\
&&\qquad\leq \frac{c{\cal C}}{\nu^\frac12}\|D(\u^1_f-\u^2_f)\|^2_{\Omega_f}.
\end{eqnarray*}

Thanks to the definition of ${\cal C}$ in Theorem \ref{th0} and (\ref{10050959}), we have
$$
c\nu\|D(\u^1_f-\u^2_f)\|^2_{\Omega_f}+c\|\mathds{K}\nabla(\phi_p^1-\phi_p^2)\|^2_{\Omega_p}\leq 0.
$$
This leads to the global uniqueness of the weak solution of ($P$) in ${\bf W}$. The uniqueness of the solution of ($Q$) in ${\bf Y}$ is obvious thanks to (\ref{LBB}). \hfill$\Box$

\section*{Conclusion} 
By introducing an auxiliary partial differential system in the porous media flow region that compensates the nonlinear convection in the Navier-Stokes equations in the energy balance, we establish an a priori estimate of weak solutions of the Navier-Stokes/Darcy model with BJS interface condition. Then an existence result of a weak solution to this coupled system is obtained without the restriction of the small data and/or the large viscosity for the first time. Finally, a global uniqueness result of the weak solution is derived which solves the open question raised in \cite{GIRAULT}.


\end{document}